\documentclass[11pt]{amsart}

\input epsf.tex
\usepackage{pst-all}
\usepackage{amsmath}
\usepackage[latin1]{inputenc}
\usepackage{amsfonts}
\usepackage{amssymb}
\usepackage{multicol}
\usepackage{verbatim}
\usepackage{amsthm}
\usepackage{amscd}
\usepackage{graphicx,epsfig}
\usepackage{pb-diagram}
\usepackage[mathscr]{eucal}
\usepackage[latin1]{inputenc}

\usepackage{amsmath}

\usepackage{amsfonts}

\usepackage{amssymb}

\usepackage{amsthm}

\usepackage{graphics}

\newcommand{\ZZ}{\mathbb{Z}}

\newcommand{\CC}{\mathbb{C}}

\newcommand{\NN}{\mathbb{N}}

\newcommand{\yl}{15pt}

\newcommand{\ffbox}[1]{
\setbox9=\hbox{$\scriptstyle\overline{1}$}
\framebox[\yl][c]{\rule{0mm}{\ht9}${\scriptstyle #1}$}
}

\newcommand{\Glie}{\mathfrak{g}}

\newcommand{\Hlie}{\mathfrak{h}}

\newcommand{\demo}{\noindent {\it \small Proof:}\quad}

\newcommand{\U}{\mathcal{U}}

\newtheorem{thm}{Theorem}[section]

\newtheorem{defi}[thm]{Definition}

\newtheorem{cor}[thm]{Corollary}

\newtheorem{prop}[thm]{Proposition}

\newtheorem{lem}[thm]{Lemma}

\newtheorem{conj}[thm]{Conjecture}

\newtheorem{rem}[thm]{Remark}

\usepackage[all]{xy}
\title{The algebra $\U_q(\hat{sl}_\infty)$ and applications}

\author{David Hernandez}

\address{CNRS -- \'Ecole Polytechnique, CMLS 91128 Palaiseau cedex - FRANCE}

\dedicatory{To Corrado De Concini on his 60th birthday}

\begin{document}

\begin{abstract} 
In this note we consider the algebra $\U_q(\hat{sl}_\infty)$ and
we study the category $\mathcal{O}$ of its integrable representations. 
The main motivations are applications to quantum toroidal algebras 
$\U_q(sl_{n+1}^{tor})$, more precisely predictions of character formulae
for representations of $\U_q(sl_{n+1}^{tor})$.
In this context, we state a general positivity conjecture for representations
of $\U_q(\hat{sl}_\infty)$ viewed as representations of $\U_q(sl_{n+1}^{tor})$,
that we prove for Kirillov-Reshetikhin modules. 
\end{abstract}

\maketitle

\vspace{.5cm}

\tableofcontents

\section{Introduction}

Consider the Dynkin diagram $X_n$ of type $A_n$.

\vspace{.5cm}

\begin{center}
\epsfig{file=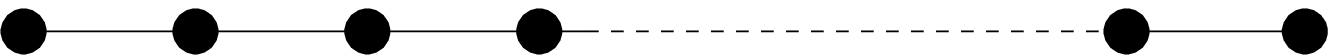,width=.5\linewidth}
\end{center}

\vspace{.5cm}

There are two limits of $X_n$ when $n\rightarrow +\infty$.
The infinite Dynkin diagram $X_{+ \infty}$ 

\vspace{.5cm}

\begin{center}
\epsfig{file=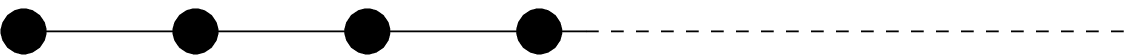,width=.5\linewidth}
\end{center}

\vspace{.5cm}

and the infinite Dynkin diagram $X_\infty$.

\vspace{.5cm}

\begin{center}
\epsfig{file=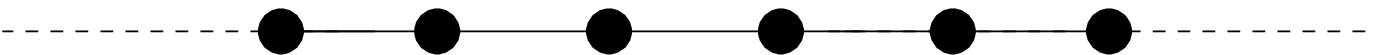,width=.5\linewidth}
\end{center}

\vspace{.5cm}

The corresponding infinite dimensional Lie algebras $sl_{+\infty}$ and $sl_{\infty}$
are defined by using the Serre presentation with an infinite root system, as in \cite{kac}.
The corresponding quantum groups $\U_q(sl_{+\infty})$ and $\U_q(sl_{\infty})$
are defined in the same way, and can be considered as limits of $\U_q(sl_{n+1})$ when $n\rightarrow \infty$. 
Such quantum groups and related structures have been studied by various authors 
(see for example \cite{ek, Fre3, jl, jmmo, ls} and references therein).

Now let us consider the analog problem for the quantum affine algebras 
$\U_q(\hat{sl}_{n+1})$. In opposition to the previous case, there is no 
obvious limit of the Dynkin diagram $X_n^{(1)}$ of type $A_n^{(1)}$ when 
$n\rightarrow \infty$.

\vspace{.5cm}

\begin{center}
\epsfig{file=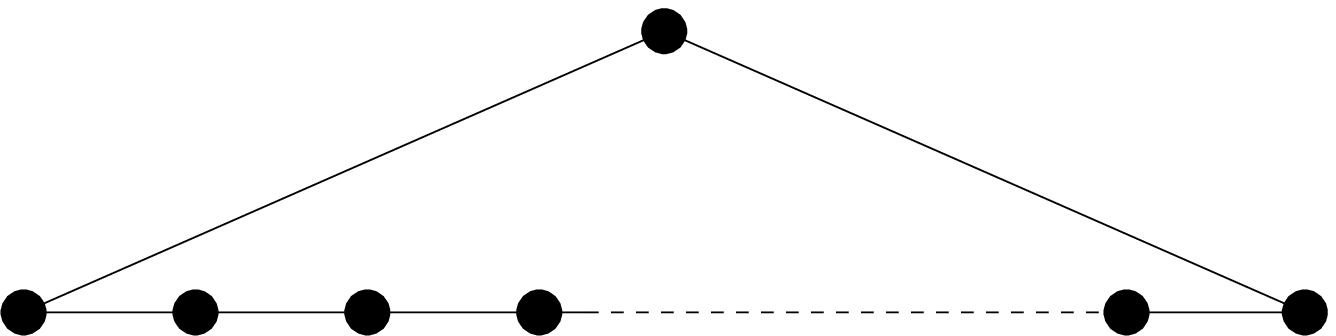,width=.5\linewidth}
\end{center}

\vspace{.5cm}

The aim of this note is threefold. First 
we aim at defining the "limit" $\U_q(\hat{sl}_\infty)$ 
of $\U_q(\hat{sl}_{n+1})$ when 
$n\rightarrow \infty$. Then, by studying the representation 
theory of $\U_q(\hat{sl}_\infty)$, we get candidates
for character formulae of representations of quantum toroidal
algebras. This is stated in a precise positivity
conjecture. Finally, we prove this conjecture for an important
classes of representations, that is Kirillov-Reshetikhin modules.

Let us explain this in more details. We first view 
$\U_q(\hat{sl}_{n+1})$ as 
a quantum affinization of $\U_q(sl_{n+1})$, that is, we use the Drinfeld 
presentation of $\U_q(\hat{sl}_{n+1})$ \cite{Dri2, bec}. Then we propose 
naturally to define the limit 
algebra $\U_q(\hat{sl}_\infty)$ as a quantum affinization of 
$\U_q(sl_{\infty})$, in the spirit of \cite{Naams, h2}.

We study the representation theory of $\U_q(\hat{sl}_\infty)$,
or more precisely of its quotient the quantum loop algebra $\U_q(\mathcal{L}sl_\infty)$.
As for quantum affinizations with a finite Dynkin diagram, the 
simple integrable modules in the category $\mathcal{O}$ are parameterized by Drinfeld polynomials,
and we can define various families of representations such as Kirillov-Reshetikhin modules.
The structure of these representations can be described by using the representation
theory of usual quantum affine algebras $\U_q(\hat{sl}_{n+1})$ 
(Theorem \ref{cond} and Proposition \ref{limit}). We derive ($q$)-character formulae for
Kirillov-Reshetikhin modules of $\U_q(\hat{sl}_\infty)$ (Theorem \ref{nexp1}).

Our main motivations are applications
to quantum toroidal algebras $\U_q(sl_{n+1}^{tor})$, that is quantum affinizations
of quantum affine algebras. Indeed, $X_\infty$ has a family of
Dynkin diagram automorphisms given by shifts.

\vspace{.5cm}

\begin{center}
\epsfig{file=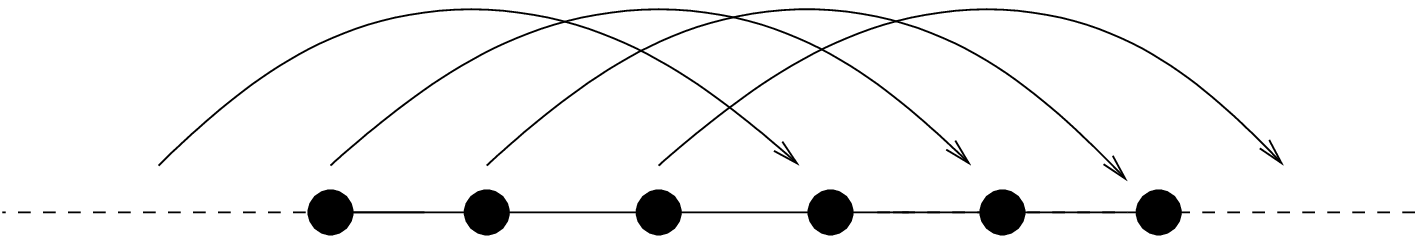,width=.5\linewidth}
\end{center}

\vspace{.5cm}

\noindent The corresponding quotients are 
the diagrams $X_n^{(1)}$. Such situations occur in the context of
twisted quantum affine algebras, whose underlying Dynkin diagram is obtained
from an automorphism of the Dynkin diagram associated to a simply-laced untwisted 
quantum affine algebra (details on this analogy will be given in the paper). 
By using arguments of \cite{h8}, this allows us to 
relate the representation theory of $\U_q(\hat{sl}_\infty)$ to the representation theory 
of $\U_q(sl_{n+1}^{tor})$, that is to predict character formulae for
representations of $\U_q(sl_{n+1}^{tor})$ from the representation theory
of $\U_q(\hat{sl}_\infty)$ (which is better understood !). 
For example we get explicit ($q$)-character formulae for
Kirillov-Reshetikhin modules of $\U_q(sl_{n+1}^{tor})$ in Theorem \ref{nexp} 
(this result was announced in \cite{Selecta}).

Let us explain this more precisely. In \cite{h5}, the author had established
a characterization of q-characters of Kirillov-Reshetikhin modules of 
$\U_q(sl_{n+1}^{tor})$. But in general, it is not clear how to extract
an explicit formula from such kind of characterization. For example, 
there is no known formula for fundamental representations of 
general quantum affinizations, although the
characterization also holds \cite{h5}. 
A point is to find candidates for character formulae.
The idea of the present paper is to study first the representation
theory of $\U_q(\hat{sl}_\infty)$ : explicit formulae
for $q$-characters of Kirillov-Reshetikhin modules are first 
established for $\U_q(\hat{sl}_\infty)$ (from quantum affine algebras
of type $A$). Then, this gives a candidate for
an explicit formula in the quantum toroidal case (by "twisting"). After having 
written the formula, we check that it satisfies the combinatorial characterization, hence
we prove the formula.

Note that it is most probably possible to prove also some results of this 
paper, especially for quantum toroidal algebras, by using geometric approach
(see references in \cite{ncluster}). But our aim is to emphasis 
the relations between $\U_q(\hat{sl}_\infty)$ and quantum toroidal algebras with 
direct methods.

In the spirit of this paper, notice that an arbitrary integrable representation
in the category $\mathcal{O}$ for $sl_\infty$ has an action
of $\hat{sl}_{n+1}$ compatible with characters (Proposition \ref{posres}; this
is also related to the study of the Fock space \cite{jmmo, vv2, tu}).
This observation in the classical case naturally leads
to conjecture an analog statement for representations of $\U_q(\hat{sl}_\infty)$.
Indeed, we prove that it can be seen as a representation of $\U_q(sl_{n+1}^{\text{tor}})$,
which is {\it a priori} virtual, and then we conjecture that it is an actual representation.
This is a positivity conjecture (Conjecture \ref{cpos}). Our results imply this conjecture for
Kirillov-Reshetikhin modules.

{\bf Acknowledgments :} The author would like to thank Edward Frenkel for its comment
in Remark \ref{knclass}.

\section{The algebra $\U_q(\hat{sl}_\infty)$}

First let us give some definitions.

We consider $C=(C_{i,j})_{i,j\in\ZZ}$ the infinite Cartan matrix associated to the infinite Dynkin diagram $X_\infty$, that is for $i,j\in \ZZ$ satisfying $i - j\notin \{-1,0,1\}$, we set 
\begin{equation}\label{cartan}C_{i,i} = 2\text{ , }C_{i,i+1} = C_{i+1,i} = -1\text{ , }C_{i,j} = 0.\end{equation}
$q\in\mathbb{C}^*$ is not a root of unity and is fixed. For $l\in\ZZ$ and $0\leq k\leq s$, we set
$$[l]_q=\frac{q^l-q^{-l}}{q-q^{-1}}\in\mathbb{Z}[q^{\pm}]\text{ , }[s]_q! = [s]_q\cdots [1]_q\text{ , }\begin{bmatrix}s\\k\end{bmatrix}_q = \frac{[s]_q!}{[s-k]_q![k]_q!}.$$

\begin{defi}\label{defiqaf} The algebra $\U_q(\hat{sl}_\infty)$ is the algebra defined by generators $x_{i,r}^{\pm}$ ($i, r\in\ZZ$), $k_i^{\pm 1}$ ($i\in\ZZ$), $h_{i,m}$ ($i\in \ZZ, m\in\ZZ\setminus\{0\}$), central elements $c^{\pm 1/2}$ and relations ($i,j\in\ZZ$, $r, r', r_1, r_2\in\ZZ$, $m, m'\in\ZZ$) :
$$[k_i, k_j] = [k_i, h_{j,m}] = 0\text{ , }[h_{i,m},h_{j,m'}] = \delta_{m, - m'} \frac{[mC_{i,j}]_q\left(c^m - c^{-m}\right)}{m\left(q - q^{-1}\right)},$$
$$k_ix_{j,r}^{\pm}k_i^{-1}=q^{\pm C_{i,j}}x_{j,r}^{\pm}
\text{ , }[h_{i,m},x_{j,r}^{\pm}] = \pm \frac{1}{m}[m C_{i,j}]_q c ^{-|m|/2} x_{j,m+r}^{\pm},$$
$$[x_{i,r}^+,x_{j,r'}^-] = \delta_{i,j}\frac{c^{(r-r')/2}\phi^+_{i,r+r'}-c^{-(r-r')/2}\phi^-_{i,r+r'}}{q-q^{-1}},$$
$$x_{i,r+1}^{\pm}x_{j,r'}^{\pm} - q^{\pm C_{i,j}}x_{j,r'}^{\pm}x_{i,r+1}^{\pm}
=q^{\pm C_{i,j}}x_{i,r}^{\pm}x_{j,r'+1}^{\pm}-x_{j,r'+1}^{\pm}x_{i,r}^{\pm},$$
$$[x_{i,r}^\pm, x_{j,r'}^\pm] = 0\text{ if $i - j\notin\{-1,0,1\}$,}$$
$$x_{i,r_1}^\pm x_{i,r_2}^\pm x_{j,r}^\pm - (q + q^{-1}) x_{i,r_1}^\pm x_{j,r}^\pm x_{i,r_2}^\pm
+ x_{j,r}^\pm x_{i,r_1}^\pm x_{i,r_2}^\pm$$
$$+ x_{i,r_2}^\pm x_{i,r_1}^\pm x_{j,r}^\pm - (q + q^{-1}) x_{i,r_2}^\pm x_{j,r}^\pm x_{i,r_1}^\pm
+ x_{j,r}^\pm x_{i,r_2}^\pm x_{i,r_1}^\pm = 0\text{ if $i-j \in\{-1,1\}$,}$$
$$\text{ where }\phi_i^\pm (z) = \underset{m\geq 0}{\sum}\phi_{i,\pm m}^{\pm}z^{\pm m} = k_i^{\pm 1}\text{exp}\left(\pm\left(q-q^{-1}\right)\underset{m'\geq 1}{\sum}h_{i,\pm m'}z^{\pm m'}\right).$$
\end{defi}

\begin{rem} These defining formulae are obtained from the general framework of quantum
affinizations (see \cite{h2} and references therein). The last two relations are analogs
of Serre relations.
\end{rem}

The subalgebra generated by the $k_i^{\pm 1}$, $x_i^{\pm} = x_{i,0}^\pm$ ($i\in\ZZ$) is isomorphic to $\U_q(sl_\infty)$.
Let $J = [a, b]\subset \ZZ$. Then the subalgebra $\hat{\U}_J$ of $\U_q(\hat{sl}_\infty)$ generated by the $x_{i,m}^{\pm}$ ($i\in J$, $m\in\ZZ$), the $k_i^{\pm 1}$ ($i\in J$), the $h_{i,r}$ ($i\in J$, $r\in\ZZ\setminus\{0\}$) and $c^{\pm 1/2}$ is isomorphic to the quantum affine algebra $\U_q(\hat{sl}_{b - a + 2})$.
Its subalgebra $\U_J$ generated by the $x_i^{\pm}$ and the $k_i^{\pm 1}$ ($i\in J$) is isomorphic to $\U_q(sl_{b - a + 2})$. The quotient of $\U_q(\hat{sl}_\infty)$ by the relation $c^{\pm 1/2} = 1$ is the quantum loop algebra $\U_q(\mathcal{L}sl_\infty)$. 

Let $\U_q(\hat{sl}_\infty)^\pm$ be the subalgebra of $\U_q(\hat{sl}_\infty)$ generated by the $x_{i,m}^{\pm}$ ($i, m\in\ZZ$) and $\U_q(\hat{\Hlie})$ the subalgebra of $\U_q(\hat{sl}_\infty)$ generated by the $k_i^{\pm 1}$ ($i\in \ZZ$), the $h_{i,r}$ ($i\in \ZZ$, $r\in\ZZ\setminus\{0\}$) and $c^{\pm 1/2}$. Then we have the following triangular decomposition.

\begin{prop}\label{dtrian} We have an isomorphism of vector spaces 
$$\U_q(\hat{sl}_\infty)\simeq \U_q(\hat{sl}_\infty)^-\otimes \U_q(\hat{\Hlie})\otimes \U_q(\hat{sl}_\infty)^+.$$
\end{prop}

\demo When $J\subset \ZZ$ is finite, an analog result is well known for the subalgebras 
$$\hat{\U}_J\simeq \hat{\U}_J^+\otimes \U_q(\hat{\Hlie}_J)\otimes \hat{\U}_J^-$$ 
with $\hat{\U}_J^\pm = \hat{\U}_J\cap \U_q(\hat{sl}_\infty)^\pm$ and $\U_q(\hat{\Hlie}_J) = \U_q(\hat{\Hlie})\cap \hat{\U}_J$ (see \cite{bec}). Consider the multiplication map 
$$\U_q(\hat{sl}_\infty)^-\otimes \U_q(\hat{\Hlie})\otimes \U_q(\hat{sl}_\infty)^+\rightarrow \U_q(\hat{sl}_\infty).$$
By using the defining relations of the algebra $\U_q(\hat{sl}_\infty)$, this map is clearly surjective. The injectivity is a direct consequence of the result for the $\hat{\U}_J$ mentioned above, as 
$$\U_q(\hat{sl}_\infty)^-\otimes \U_q(\hat{\Hlie})\otimes \U_q(\hat{sl}_\infty)^+ = \bigcup_{J\subset \ZZ, \text{ J finite}} \hat{\U}_J^+\otimes \U_q(\hat{\Hlie}_J)\otimes \hat{\U}_J^-.$$
\qed

\begin{rem}\label{hcom} The image $\U_q(\mathcal{L}\Hlie)$ of $\U_q(\hat{\Hlie})$ in $\U_q(\mathcal{L}sl_\infty)$ is commutative.
\end{rem}

\section{Representations of $\U_q(\hat{sl}_\infty)$}

We remind standard results on the representation theory of $\U_q(sl_\infty)\subset \U_q(\hat{sl}_\infty)$.
Then we study representations of $\U_q(\hat{sl}_\infty)$ by using the representation theory of quantum affine algebras $\U_q(\hat{sl}_{n+1})$. We deduce a character formula for a class of representations called Kirillov-Reshetikhin modules.

\subsection{Representations of $\U_q(sl_\infty)$}\label{intrep} First let us remind standard definitions (for instance, see \cite{kac, cp1} and references therein).

Let $P = \bigoplus_{i\in\ZZ} \ZZ \Lambda_i$ where the $\Lambda_i$ are the fundamental weights. For $i\in\ZZ$, we have the simple root 
$$\alpha_i = 2\Lambda_i - \Lambda_{i+1} - \Lambda_{i-1}\in P.$$ 
Let $P^+ = \bigoplus_{i\in\ZZ}\NN \Lambda_i\subset P$ and $Q^+={\sum}_{i\in I}\mathbb{N} \alpha_i$.

First let us remind known facts about the representation theory of $\U_q(sl_{\infty})$, which are very similar to results for $\U_q(sl_{n+1})$. For $V$ a representation of $\U_q(sl_{\infty})$ and $\nu = \sum_{i\in I}\nu_i \Lambda_i \in P$, the weight space $V_{\nu}$ is :
$$V_\nu = \{v\in V|k_i.v = q^{\nu_i}v,\forall i\in \ZZ\}.$$
We say that $V$ is $\Hlie$-diagonalizable if $V = \bigoplus_{\nu\in P} V_\nu$ and $V_\nu$ is finite-dimensional for any $\nu\in P$. For such a representation, we can define its character as the formal infinite sum
$$\chi(V) = \sum_{\lambda\in P} \text{dim}(V_\lambda) e(\lambda),$$
where the $e(\lambda)$ are formal elements depending on $\lambda\in P$. 
Let us define the category $\mathcal{O}$ as in \cite{kac}.

\begin{defi} $V$ is said to be in the category $\mathcal{O}$ if $V$ is $\Hlie$-diagonalizable and if there are $\lambda_1,\cdots,\lambda_N\in P$ such that
$$\{\nu\in P|V_\nu\neq\{0\}\}\subset \bigcup_{j=1\cdots N} (\lambda_j - Q^+).$$
\end{defi}

\noindent A representation $V$ is said to be of highest weight $\lambda\in P$ if there is $v\in V_\lambda$, called an highest weight vector, satisfying $\U_q(sl_\infty).v = V$ and $x_i^+.v = 0$ for any $i\in\ZZ$. 

\begin{lem}\label{unionun} For $V$ of highest weight $\lambda$ with highest vector $v$, we have 
$$V = \bigcup_{n\geq 0} (\U_{[-n,n]}v).$$
As a consequence $V$ is in the category $\mathcal{O}$. Moreover,
if $V$ is simple, then for $n\geq 0$, $\U_{[-n,n]}v$ is simple as a $\U_{[-n,n]}$-module.
\end{lem}

\demo The first point follows from $\U_q(sl_\infty) = \bigcup_{n\geq 0}\U_{[-n,n]}$. Each 
$\U_{[-n,n]}v$ is clearly of highest weight. 
From the well-known structure of highest weight modules of $\U_{[-n,n]}\simeq \U_q(sl_{2(n + 1)})$, $V$ is $\Hlie$-diagonalizable. Moreover, $V = \U_q(sl_\infty)^-.v$ where $\U_q(sl_\infty)^-$ is the subalgebra of $\U_q(sl_\infty)$ generated by the $x_i^-$ ($i\in\ZZ$). Hence a weight $\nu\in P$ such that $V_\nu\neq \{0\}$ satisfies $\nu\in \lambda - Q^+$. So $V$ is in the category $\mathcal{O}$. 

Now suppose that $V$ is simple. Then there is no highest weight vector
in $V\setminus \CC v$. For $w\in \left(\U_{[-n,n]}v\right)\setminus \CC v$, we have $\sum_{i\notin [-n,n]}\CC x_i^+w = 0$ for weight reason. So if $w$ is of highest weight in $\U_{[-n,n]}v$, then $w$ is of highest weight in $V$, contradiction. So $\U_{[-n,n]}v$ is simple.\qed

By standard highest weight theory arguments, for each $\lambda\in P$, there is a unique simple highest weight module $V(\lambda)$ of highest weight $\lambda$. These are the simple objects of the category $\mathcal{O}$.

\begin{defi} $V$ is said to be integrable if $V$ is $\Hlie$-diagonalizable and if 
$$V_{\nu \pm N\alpha_i} = \{0\}\text{ for all $\nu\in P$, $N>>0$, $i\in \ZZ$.}$$
\end{defi}

\begin{thm}\label{fincase}
The simple integrable representations of $\U_q(sl_\infty)$ in the category $\mathcal{O}$ are the $V(\lambda)$ where $\lambda\in P^+$.
\end{thm}

This is a well-known result, we give a proof for completeness (and as analog arguments will be used later).

\demo First such a representation $V(\lambda)$ is in the category $\mathcal{O}$ by Lemma \ref{unionun}. 
Let $v$ be a highest weight vector of $V(\lambda)$. For $\nu \in P$ and $i\in\ZZ$, by Lemma \ref{unionun} there is $n\geq 0$ such that $|i| < n$ and $V_\nu \subset \U_{[-n,n]}v$. Moreover $\U_{[-n,n]}v$ is a simple $\U_{[-n,n]}$-module by Lemma \ref{unionun}. Then for any $N\geq 0$, $V_{\nu \pm N\alpha_i}\subset \U_{[-n,n]}v$. Hence, from the result for $\U_{[-n,n]}$ (see for example \cite{lu}), we get $V_{\nu \pm N\alpha_i} = \{0\}$ for $N>>0$. So $V$ is integrable.

Now for $V$ a simple integrable representation in the category $\mathcal{O}$, there is $\lambda = \sum_{i\in\ZZ} n_i \Lambda_i \in P$ such that $V \simeq V(\lambda)$. Then
for each $i\in\ZZ$, we have $n_i\geq 0$ from the analog result for $\U_{\{i\}}\simeq \U_q(sl_2)$ ($n_i\Lambda_i$ is the highest weight of an integrable $\U_q(sl_2)$-module). So $\lambda\in P^+$.
\qed

Note that, in opposition to the case of $\U_q(sl_{n+1})$, the algebra $\U_q(sl_{\infty})$ has simple integrable representations which are not highest of lowest weight. For example, consider the following representation, inspired by Kashiwara extremal representations \cite{kasm1}. We define the action of $\U_q(sl_{\infty})$ on $V = \bigoplus_{i\in \ZZ} \CC v_i$ by the following formulae ($i,j\in\ZZ$) :
$$k_i . v_j = q^{\delta_{i,j}-\delta_{i,j+1}} v_j\text{ , }x_i^+.v_j = \delta_{i,j+1} v_i\text{ , }x_i^-.v_j = \delta_{i,j}v_{i+1}.$$
This representation is clearly simple and integrable, and
$$\chi(V) = \sum_{i\in\ZZ}e(\Lambda_i - \Lambda_{i-1}).$$
In particular $V$ is not in the category $\mathcal{O}$.

\subsection{Representations of $\U_q(\hat{sl}_\infty)$}

We give various definitions and results mimicking \cite{m2, Naams, h2} for quantum affinizations. 
We focus on representations with a trivial action of 
$c^{\pm \frac{1}{2}}$, that is representations of $\U_q(\mathcal{L}sl_\infty)$.

\begin{defi} A representation of $\U_q(\hat{sl}_\infty)$ is said to be integrable (resp. in the category $\mathcal{O}$) if it is integrable (resp. in the category $\mathcal{O}$) as a $\U_q(sl_\infty)$-module.\end{defi}

Let $\mathcal{O}_{\text{int}}$ be the category of integrable representations of $\U_q(\hat{sl}_\infty)$ in the category $\mathcal{O}$. Let $\text{Spec}(\U_q(\mathcal{L}\Hlie))$ be the set of algebra morphisms 
$$\gamma : \U_q(\hat{\Hlie})\rightarrow \CC$$ 
satisfying $\gamma(c^{\pm\frac{1}{2}}) = 1$ and $\gamma(k_i)\in q^\ZZ$ for any $i\in\ZZ$. 

For $\gamma\in \text{Spec}(\U_q(\mathcal{L}\Hlie))$, we define $\omega(\gamma) = \sum_{i\in I}n_i\Lambda_i \in P$ where the $n_i\in\ZZ$ are set so that $\gamma(k_i) = q^{n_i}$.    
    
\begin{defi}A representation $V$ of $\U_q(\hat{sl}_\infty)$ is said to be of $l$-highest weight $\gamma\in \text{Spec}(\U_q(\mathcal{L}\Hlie))$ if there is $v\in V\setminus\{0\}$ (called $l$-highest weight vector) such that 

i) $V = \U_q(\hat{sl}_\infty)^-.v$, 

ii) for any $h\in \U_q(\hat{\Hlie})$, $h.v = \gamma(h) v$, 

iii) for any $i\in I, m\in\ZZ$, $x_{i,m}^+.v=0$.
\end{defi}

\noindent For any $\gamma\in \text{Spec}(\U_q(\mathcal{L}\Hlie))$, we have by Proposition \ref{dtrian} a corresponding Verma module $M(\gamma)$ of $l$-highest weight $\gamma$. By standard arguments, it has a unique simple quotient $L(\gamma)$.

\begin{lem}\label{unionbis} Let $v_\gamma$ be a $l$-highest weight vector of $L(\gamma)$. Then 
$$L(\gamma) = \bigcup_{n\geq 0} L_n(\gamma)\text{ where }L_n(\gamma)=\hat{\U}_{[-n,n]} v_\gamma.$$
\end{lem}

\demo As in Lemma \ref{unionun}, we get the result from
$$\U_q(\hat{sl}_\infty) = \bigcup_{n\geq 0} \hat{\U}_{[-n,n]}.$$ 
\qed

In opposition to the case $\U_q(sl_\infty)$ in Lemma \ref{unionun}, the representation
$L(\gamma)$ is not necessarily in the category $\mathcal{O}$ (see \cite[Lemma 4.11]{h2}).
However, we have the following.

\begin{thm}\label{cond} The simple objects of $\mathcal{O}_{\text{int}}$ are the $L(\gamma)$ such that there is $(P_i)_{i\in\ZZ}\in (1+u\mathbb{C}[u])^{(\ZZ)}$ satisfying for $i\in \ZZ$ the relation in $\mathbb{C}[[z]]$ (resp. in $\mathbb{C}[[z^{-1}]]$):
$$\gamma(\phi_i^\pm(z))=q^{\pm \text{deg}(P_i)}\frac{P_i(zq^{-1})}{P_i(zq)}.$$
Each $L_n(\gamma)$ is a simple finite dimensional $\hat{\U}_{[-n,n]}$-module.
\end{thm}
Here the notation $(1+u\mathbb{C}[u])^{(\ZZ)}$ means that only a finite number of $P_i$ are not equal to $1$.
The $P_i$ are called Drinfeld polynomials of $L(\gamma)$.

\demo The fact that the $L_n(\gamma)$ are simple is proved as in \cite[Lemma 5.9]{AnnENS} : for $v\in L_n(\gamma)$, we have $\sum_{i\notin [-n,n], m\in\ZZ}x_{i,m}^+v = 0$ for weight reason. So if it is not of highest weight in $L(\gamma)$, there is $i\in [-n,n], m\in\ZZ$ so that $x_{i,m}^+v\neq 0$, hence $v$ is not of highest weight in $L_n(\gamma)$. 
Now from the form of $\gamma$ given in Theorem \ref{cond}, $L_n(\gamma)$ is finite dimensional by \cite[Theorem 12.2.6]{cp1}. Now we can check as in Lemma \ref{unionun} that $L(\gamma)$ is in the category $\mathcal{O}$, and we conclude as in the proof of Theorem \ref{fincase} by using Lemma \ref{unionbis}.\qed

Example :  Let $k\geq 0$, $a\in\CC^*$, $i\in\ZZ$, and consider $\gamma_{k,a}^{(i)}\in\text{Spec}(\U_q(\mathcal{L}\Hlie))$ associated to the Drinfeld polynomials
\begin{equation*}
\begin{split}
P_j(u) =
\begin{cases}
(1-ua)(1-uaq^2)\cdots(1-uaq^{2(k - 1)})&\text{ for $j = i$,}
\\1&\text{ for $j\neq i$.}
\end{cases}
\end{split}
\end{equation*}
$W_{k,a}^{(i)} = L(\gamma_{k,a}^{(i)})$ is called a Kirillov-Reshetikhin module.

\begin{rem}\label{eval} Note that for $n\geq 0$, $L_n(\gamma_{k,a}^{(i)})$ 
is a usual Kirillov-Reshetikhin module of $\U_q(\hat{sl}_{2(n+1)})$, that is
it can be constructed from the simple $\U_q(sl_{2(n+1)})$-module of highest $k\Lambda_i$ by using a Jimbo
evaluation morphism $\U_q(\hat{sl}_{2(n+1)})\rightarrow \U_q(sl_{2(n+1)})$ 
(see for instance \cite{cp1}). In particular each $L_n(\gamma_{k,a}^{(i)})$
is simple as a representation of $\U_{[-n,n]}\simeq \U_q(sl_{2(n+1)})$, hence $W_{k,a}^{(i)}$ is simple
as a $\U_q(sl_\infty)$-module (although a priori it is not clear, at least for the author,
how to define an evaluation morphism $\U_q(\hat{sl}_\infty)\rightarrow \U_q(sl_\infty)$).
\end{rem}

\subsection{Characters}

Let $V$ in $\mathcal{O}_{\text{int}}$. From Remark \ref{hcom}, we have 
$$V = \bigoplus_{\gamma\in \text{Spec}(\U_q(\mathcal{L}\Hlie))}V_\gamma,$$ 
$$\text{ where }V_\gamma = \{v\in V|\exists p \geq 0, \forall h\in\U_q(\hat{\Hlie}), (h - \gamma(h))^p.v = 0\}.$$
We have $V_\gamma\subset V_{\omega(\gamma)}$, that is we get a finer decomposition than the decomposition in weight spaces.
We define the $q$-character of $V$ as in \cite{Fre} 
$$\chi_q(V) = \sum_{\gamma\in \text{Spec}(\U_q(\mathcal{L}\Hlie))} \text{dim}(V_\gamma) e(\gamma),$$
where $e(\gamma)$ is a formal element. The formal sum $\chi_q(V)$ belongs to $\mathcal{E}$ which is defined \cite{h2} as the group of formal infinite sums $\sum_\gamma n_\gamma e(\gamma)$ such that 

1) for any $\lambda\in P$ the set $\{\gamma|n_\gamma\neq 0\text{ , }\omega(\gamma) = \lambda\}$ is finite,

2) there are $\lambda_1,\cdots, \lambda_p\in P$ such that 
$$\{\omega(\gamma)|\gamma, n_\gamma\neq 0\}\subset \bigcup_{j=1\cdots p} (\lambda_j - Q^+).$$

\noindent We get an injective group morphism
$$\chi_q :\text{Rep}(\U_q(\hat{sl}_\infty))\rightarrow \mathcal{E}$$
where $\text{Rep}(\U_q(\hat{sl}_\infty))$ is the Grothendieck group of $\mathcal{O}_{\text{int}}$.

\begin{rem}\label{sdtchar} If we replace \cite{Fre} each $e(\gamma)$ by $\omega(\gamma)$, we get $\chi(V)$ from $\chi_q(V)$.
\end{rem}

For $J\subset \ZZ$ finite and $V$ a finite dimensional representation of $\hat{\U}_J$, we have $\chi_q(V)\in\bigoplus_{\gamma : \hat{\U}_J\cap \U_q(\hat{\Hlie})\rightarrow \CC}\ZZ e(\gamma)$ its usual $q$-character.

Let $V$ in $\mathcal{O}_{\text{int}}$. Let $n\geq 0$. Let us write the (usual) $q$-character of $V_{[-n,n]}$ as a $\hat{\U}_{[-n,n]}$-module
$$\chi_q(V_{[-n,n]}) = \sum_{\gamma, \omega(\gamma) = \lambda_n} N_\gamma e(\gamma).$$
Each $\gamma : \hat{\U}_{[-n,n]}\cap \U_q(\hat{\Hlie})\rightarrow \CC$ can be extended to a unique algebra morphism $\tilde{\gamma}\in \text{Spec}(\U_q(\mathcal{L}\Hlie))$ by setting $\tilde{\gamma}(h_{i,r}) = 0$, $\tilde{\gamma}(k_i) = 1$ for $|i| > n$, $r\neq 0$. Let $\lambda\in P$ and consider the image $\lambda_n$ of $\lambda$ in the weight lattice of $\hat{\U}_{[-n , n]}$. We get an element
$$\chi_{q,\lambda,n}(V) = \sum_{\gamma, \omega(\gamma) = \lambda_n} N_\gamma e(\tilde{\gamma})\in \mathcal{E}.$$

\begin{prop}\label{limit} For each $\lambda\in P$, the sequence $(\chi_{q,\lambda,n}(V))_{n\geq 0}\in \mathcal{E}^\NN$ is stationary when $n\rightarrow \infty$. Let $\chi_{q,\lambda}(V)$ be its limit. Then
$$\chi_q(V) = \sum_{\lambda\in P}\chi_{q,\lambda}(V).$$
\end{prop}

\demo Let $\lambda\in P$. From Lemma \ref{unionbis}, there is $n\geq 0$ such that $V_\lambda\subset V_{[-n,n]}$. Moreover there is $n'\geq n$ such that $\lambda, \mu \in \sum_{-n'\leq j\leq n'}\ZZ \Lambda_j$ where $\mu$ is the highest weight of $V$. Now for any $N\in\ZZ$ such that $|N|> n'$, we have
$$[h_{N,r},\hat{\U}_{[-n,n]}] = [k_N,\hat{\U}_{[-n,n]}] = 0$$
for any $r\neq 0$. Besides the action of $h_{N,r}$ on $V_\mu$ is $0$ and the action of $k_N$ is the identity. So this is the same on $V_\lambda$. This implies that
$$\chi_{q,\lambda,N}(V) = \chi_{q,\lambda}(V) = \sum_{\{\gamma\in \text{Spec}(\U_q(\mathcal{L}\Hlie))|V_\gamma\subset V_\lambda\}} \text{dim}(V_\gamma)e(\gamma).$$
\qed

\begin{cor} 
Let $V$ in $\mathcal{O}_{\text{int}}$ and $\gamma\in\text{Spec}(\U_q(\mathcal{L}\Hlie))$ such that $V_\gamma\neq 0$.
For any $i\in \ZZ$, there are $Q_i(z),R_i(z)\in\mathbb{C}[z]$ with constant term equal to $1$ such that
$$\gamma(\phi_i^\pm (z)) = q^{\text{deg}(Q_i)-\text{deg}(R_i)}\frac{Q_i(zq^{-1})R_i(zq)}{Q_i(zq)R_i(zq^{-1})}.$$
\end{cor}

\demo The result is known \cite{Fre} for the quantum affine algebras $\U_q(\hat{sl}_{n+1})$. Hence the result follows from Proposition \ref{limit}.\qed

\noindent So, as in \cite{Fre}, we can define the monomial $m_\gamma = \prod_{i\in I, a\in\CC^*}Y_{i,a}^{\mu_{i,a} - \nu_{i,a}}$ where
$$Q_i(z)={\prod}_{a\in\mathbb{C}^*}(1-za)^{\mu_{i,a}}\text{ and }R_i(z)={\prod}_{a\in\mathbb{C}^*}(1-za)^{\nu_{i,a}}.$$
$e(\gamma)$ is also denoted by $m_\gamma$.

For $i\in \ZZ$ let
$$\mathfrak{K}_i = \ZZ[Y_{i,a}(1 + A_{i,aq}^{-1}),Y_{j,a}^{\pm 1}]_{a\in\CC^*,j\neq i}\text{ where }A_{i,a} = Y_{i,aq^{-1}}Y_{i,aq}Y_{i+1,a}^{-1}Y_{i-1,a}^{-1}.$$

\begin{cor}\label{whatim} For $V$ in $\mathcal{O}_{\text{int}}$ and $i\in\ZZ$, $\chi_q(V)\in \mathcal{E}$ is an infinite sum of elements in $\mathfrak{K}_i$.
\end{cor}

\demo The result holds for the subalgebras $\hat{\U}_J$ when $J$ is finite by \cite{Fre, Fre2}, and so the result follows from Proposition \ref{limit}.
\qed

Note that there are simple integrable representations of $\U_q(\hat{sl}_\infty)$ not in the category $\mathcal{O}$ with $q$-character satisfying the condition of Corollary \ref{whatim}. For example, let $V = \bigoplus_{i\in \ZZ} \CC v_i$. We define the action of $\U_q(\hat{sl}_\infty)$ by
($r\in\mathbb{Z}$, $m > 0$, $i,j \in\mathbb{Z}$) :
$$x_{i,r}^+.v_{i+j} = \delta_{j,1}q^{ri}v_i\text{ , }x_{i,r}^-.v_{i+j} = \delta_{j,0} q^{ri}v_{i+1},$$
$$\phi_{i,\pm m}^{\pm}.v_{i+j} = \pm (\delta_{j,0} - \delta_{j,1}) (q - q^{-1})q^{\pm mi }v_{i + j},$$
$$k_i^{\pm}.v_{i+j} = q^{\pm(\delta_{j,0} - \delta_{j,1})}v_{i+j},$$
This representation is clearly simple and integrable, and its $q$-character is
$$\chi_q(V) = \sum_{i\in\ZZ}\ffbox{i}_1,$$
where for $a\in\CC^*$ and $\alpha\in \ZZ$ we set $\ffbox{\alpha}_a = Y_{\alpha - 1,aq^\alpha}^{-1} Y_{\alpha,aq^{\alpha - 1}}$. 

\subsection{A character formula}

For $I,J\subset\ZZ$, a tableaux $(T_{\alpha,\beta})_{\alpha\in I, \beta\in J}$ is said to be semi-standard if it has
coefficients in $\ZZ$ which increase relatively to $\beta$ and strictly increase relatively to $\alpha$.

Consider $i\in\ZZ$ and $k\geq 0$. Let $\mathcal{T}_{i,k}$ be the set of semi-standard tableaux $(T_{\alpha,\beta})_{\alpha\leq i, 1\leq \beta\leq k}$ such that for any $1\leq \beta\leq k$, $T_{\alpha,\beta} = \alpha$ for $\alpha << 0$.

\begin{thm}\label{nexp1} We have the following explicit formula 
$$\chi_q(W_{k,a}^{(i)}) = \sum_{T\in \mathcal{T}_{i,k}} m_T$$
$$\text{ where }m_T = \prod_{\alpha\leq i,1\leq \beta\leq k} \ffbox{T_{\alpha,\beta}}_{aq^{i - 1 + 2(\beta - \alpha)}}.$$
\end{thm}

\demo Explicit $q$-character formulae are known for Kirillov-Reshetikhin modules of $\U_q(\hat{sl}_{n+1})$,
see the various references in the introduction of \cite{h4} (for instance \cite{Fre3}). We get
$$\chi_q((W_{k,a}^{(i)})_{[i-n,i+n]}) = \sum_{T\in \mathcal{T}_{i,k}^{(n)}} \prod_{\alpha\leq i,1\leq \beta\leq k} \ffbox{T_{\alpha,\beta}}_{aq^{i - 1 + 2(\beta - \alpha)}},$$
where $\mathcal{T}_{i,k}^{(n)}$ is the set of semi-standard tableaux $(T_{\alpha,\beta})_{i-n\leq \alpha\leq i, 1\leq \beta\leq k}$ with coefficients in $[i-n , i+n+1]$. In this formula, for $b\in\CC^*$ we have set $Y_{i-n-1,b} = Y_{i+n+1,b} = 1$. 

Now the result follows from Proposition \ref{limit}.
\qed

Note that the highest monomial of the right terms is obtained for $T_0 = (\alpha)_{\alpha,\beta}$. Indeed, we get
$$m_{T_0} 
= \prod_{1\leq \beta\leq k} \prod_{\alpha\leq i}Y_{\alpha - 1,aq^{-\alpha+2\beta + i -1}}^{-1} Y_{\alpha,aq^{-\alpha + 2\beta + i - 2}}= \prod_{1\leq \beta\leq k} Y_{i,aq^{2(\beta - 1)}},$$
which is by definition the highest monomial of $W_{k,a}^{(i)}$.

\section{Application to quantum toroidal algebras}\label{torapp}

The quantum toroidal algebra $\U_q(sl_{n+1}^{tor})$ ($n\geq 1$) were introduced in \cite{gkv}
(see also \cite{tu, vv2} and the references in the review paper \cite{Selecta}). One
of the main motivation for their study are connections to double affine Hecke algebras \cite{che3}.
The aim of this section is to explain how the representation theory of $\U_q(\hat{sl}_\infty)$
can be used to predict ($q$)-character formulae for representation of $\U_q(sl_{n+1}^{tor})$.
This is, with Conjecture \ref{cpos}, one of our main applications 
to the representation theory of $\U_q(sl_{n+1}^{tor})$.

Let us recall the definition of quantum toroidal algebras. Let $I_n = \ZZ/(n+1)\ZZ$.

First suppose that $n\geq 2$. The matrix $C = (C_{i,j})_{i,j\in I_n}$ is defined as above by formula (\ref{cartan}) with $i,j\in I_n$. The algebra $\U_q(sl_{n+1}^{tor})$ is defined by the same generators and relations as in Definition \ref{defiqaf} with $i,j\in I_n$.

The quantum toroidal algebra $\U_q(sl_2^{tor})$ requires a special definition. 
Consider the algebra $\tilde{\U}_q(sl_2^{tor})$ defined by generators $x_{0,r}^{\pm}, x_{1,r}^{\pm}$ ($r\in\ZZ$), $k_0^{\pm 1}, k_1^{\pm 1}$, $h_{1,m}, h_{0,m}$ ($m\in\ZZ-\{0\}$), central elements $c^{\pm 1/2}$ and relations
$$[k_i, k_j] = [k_i, h_{j,m}] = 0\text{ , }[h_{i,m},h_{i,m'}] = \delta_{m, - m'} \frac{[2m]_q}{m} \frac{c^m - c^{-m}}{q - q^{-1}},$$
$$[h_{i,m},h_{j,m'}] = \delta_{m, - m'} \frac{-2[m]_q}{m} \frac{c^m - c^{-m}}{q - q^{-1}}\text{ if $i\neq j$,}$$
$$k_ix_{j,r}^{\pm}k_i^{-1}=q^{\pm C_{i,j}}x_{j,r}^{\pm}\text{ , }[h_{i,m},x_{i,r}^{\pm}] = \pm \frac{1}{m}[2m]_q c ^{-|m|/2} x_{i,m+r}^{\pm},$$
$$[h_{i,m},x_{j,r}^{\pm}] = \mp \frac{2}{m}[m]_q c ^{-|m|/2} x_{j,m+r}^{\pm}\text{ if $i\neq j$,}$$
$$[x_{i,r}^+,x_{j,r'}^-] = \delta_{i,j}\frac{c^{(r-r')/2}\phi^+_{i,r+r'}-c^{-(r-r')/2}\phi^-_{i,r+r'}}{q-q^{-1}},$$
$$\text{ where }\phi_i^\pm (z) = \underset{m\geq 0}{\sum}\phi_{i,\pm m}^{\pm}z^{\pm m} = k_i^{\pm 1}\text{exp}(\pm(q-q^{-1})\underset{m'\geq 1}{\sum}h_{i,\pm m'}z^{\pm m'}).$$
Then $\U_q(sl_2^{tor})$ is a quotient of $\tilde{\U}_q(sl_2^{tor})$ : as for other quantum toroidal algebras, there are additional relations (analog to Serre relations) involving generators $x_{i,r}^+$ (resp. $x_{i,r}^-$). We do not write them because, as for usual highest weight theory, the simple integrable representations in the category $\mathcal{O}$ are the same for $\tilde{\U}_q(sl_2^{tor})$ and $\U_q(sl_2^{tor})$.

\begin{rem} We use here the definition of $\U_q(sl_2^{tor})$ as in \cite{ncluster}, that is we use the quantized Cartan matrix
\begin{equation}\label{cq}C(q) = \begin{pmatrix}q + q^{-1} & -2\\ -2 & q + q^{-1}\end{pmatrix}.\end{equation} 
This is different than other possible quantized Cartan matrices such as
$$\begin{pmatrix}q + q^{-1} & -(q + q^{-1})\\ -(q + q^{-1}) & q + q^{-1}\end{pmatrix}\text{ , }
\begin{pmatrix}q^2 + q^{-2} & -(q + q^{-1})\\ -(q + q^{-1}) & q^2 + q^{-2}\end{pmatrix}.$$
The first one is obtained by extrapolating the formulae in \cite{Fre} and is used for example in \cite{Naams}, but is not invertible, and the second one is used for example in \cite{h5}. The matrix $C(q)$ in formula (\ref{cq}) is introduced in \cite[Remark 3.13]{ncluster}. As it allows to state uniformly the results of this note, this is another indication that $C(q)$ is natural to define $\U_q(sl_2^{tor})$.\end{rem}

Let us go back to the general quantum toroidal algebras $\U_q(sl_{n+1}^{tor})$. The subalgebra $\U^h_q$ of $\U_q(sl_{n+1}^{tor})$ generated by the
$x_{i,0}^{\pm}$, $k_i^{\pm 1}$ for $i\in I_n$ is isomorphic to the quantum affine algebra
$\U_q(\hat{sl}_{n+1})$.

We can define as above $q$-characters depending on variables
$Y_{i,a}^{\pm 1}$ ($i\in I_n, a\in\CC^*$) for integrable representations in the category $\mathcal{O}$
with value in a group $\mathcal{E}_n$ (see \cite{Selecta} for references and details when $n\geq 2$; the definition are the same when $n = 1$). The Kirillov-Reshetikhin modules $\mathcal{W}_{k,a}^{(i)}$ are also defined in an analog way
(it should not be confused with the Kirillov-Reshetikhin module $W_{k,a}^{(i)}$ for $\U_q(\hat{sl}_\infty)$).

Consider the ring morphism 
$$\phi_n : \ZZ[Y_{i,a}^{\pm 1}]_{i\in \ZZ, a\in\CC^*}\rightarrow \ZZ[Y_{i,a}^{\pm 1}]_{i\in I_n, a\in\CC^*}$$ 
defined, for $i\in \ZZ$, $a\in\CC^*$, by  
$$\phi_n(Y_{i,a}) = Y_{[i],a}$$ 
where $[i]$ is the image of $i\in\ZZ$ in $I_n$. This morphism gives rise naturally to a ring morphism $\phi_n : \text{Im}(\chi_q)\rightarrow \mathcal{E}_n$.

We define the $A_{r,a}$ for $\U_q(sl_{n+1}^{tor})$ with the same formula as for $\U_q(\hat{sl}_\infty)$ when $n\geq 2$. For $\U_q(sl_2^{tor})$, we set $A_{r,a} = Y_{r,aq^{-1}}Y_{r,aq}Y_{r+1,a}^{-2}$.

A monomial $m$ is said to be dominant if it involves only positive powers of the variables, that is $m\in\ZZ[Y_{i,a}]_{i\in I_n,a\in\CC^*}$.

\begin{thm}\label{nexp} Let $i\in I_n$, $a\in\CC^*$. For $j\in\ZZ$ such that $[j] = i$, we have
$$\chi_q(\mathcal{W}_{k,a}^{(i)}) = \phi_n(\chi_q(W_{k,a}^{(i)}))= \sum_{T\in \mathcal{T}_{j,k}} \phi_n(m_T).$$
\end{thm}

The two right terms are equal by Theorem \ref{nexp1}. So the $q$-character formula
for $\mathcal{W}_{k,a}^{(i)}$ is {\it a priori} conjectured from the representation of $\U_q(\hat{sl}_\infty)$.
It is a first example
of the relation between the two representation theories. Note that it also provides a relation between the representation theory
of the quantum toroidal algebras $\U_q(sl_{n+1}^{tor})$ for various $n$ (as, in the formula, only $\phi_n$ depends
on $n$).

\demo From \cite[Theorem 6.14]{h5}, $\chi_q(\mathcal{W}_{k,a}^{(i)})$ is characterized by the following properties

i) it is an infinite sum of elements in $\ZZ[(Y_{r,a}(1 + A_{r,aq}^{-1}), Y_{r',a}^{\pm 1}]_{a\in\CC^*, r'\neq r}$ for each $r\in I_n$,

ii) the highest monomial $m_{T_0}$ is its unique dominant monomial.

\noindent (This is proved in \cite{h5} when $n\geq 2$, but the proof is the same for $\U_q(sl_2^{tor})$).
Hence it suffices to prove that $\phi_n(\chi_q(W_{k,a}^{(j)}))$ satisfies these properties.

For the first point, let $r\in I_n$. We have seen in Corollary \ref{whatim} that for each $R\in\ZZ$ such that $[R] = r$, $\chi_q(W_{k,a}^{(i)})$ is a sum of elements in
$$\ZZ[(Y_{R,a}(1+A_{R,aq}^{-1}), Y_{R',a}^{\pm 1}]_{a\in\CC^*, R'\neq R}.$$ 
So $\chi_q(W_{k,a}^{(i)})$ is a sum of elements in 
$$\ZZ[(Y_{R,a}(1+A_{R,aq}^{-1}), Y_{R',a}^{\pm 1}]_{a\in\CC^*, R,R'\in\ZZ, [R] = i, [R'] \neq i}.$$ 
Hence the result.

For the second point, consider $T\in \mathcal{T}_{j,k}\setminus \{T_0\}$. We have 
$$m_T = m_{T_0}A_{i,aq^{2k-1}}^{-1}B$$
where $B$ is a product of various $A_{i,q^{s-1}}^{-1}$
where $s\in\ZZ$. Let $S$ be the maximal of all such $s$ and of $2k$. Then the $Y_{j,s}$ occurring in $m_T$
have negative powers, and so $m_T$ is not dominant. For the same reason, $\phi_n(m_T)$ is not dominant.
\qed

\begin{rem} For $k = 1$, explicit formulae have been proved in \cite{nag}.\end{rem}

\noindent Let us explain the analogy of this result with the study of twisted quantum affine algebras $\U_q(\hat{\Glie}^\sigma)$ in \cite{h8}. The algebra $\U_q(\hat{\Glie}^\sigma)$
is associated to a simply-laced quantum affine algebra $\U_q(\hat{\Glie})$ and to a non trivial automorphism $\sigma$
of the underlying finite type Dynkin diagram. One of the main results of \cite{h8} is that the twisted $q$-character of a Kirillov-Reshetikhin module of $\U_q(\hat{\Glie}^\sigma)$ can be obtained from the $q$-character of a Kirillov-Reshetikhin module of $\U_q(\hat{\Glie})$ by using a certain projection of the variables $Y_{i,a}$.
The situation studied in this paper is analogous, with $\U_q(\hat{sl}_\infty)$ playing the role of $\U_q(\hat{\Glie})$ and $\U_q(sl_{n+1}^{tor})$ playing the role of $\U_q(\hat{\Glie}^\sigma)$ : we have a non trivial automorphism $\tau_n : i\mapsto i+n+1$ of the Dynkin diagram $X_\infty$. The map $\phi_n$ corresponds to the projection of variables in this case.

\begin{cor} $\mathcal{W}_{k,a}^{(i)}$ is not simple as a $\U_q^h$-module.
\end{cor}

\demo First suppose that $n$ is odd. Let $M$ be the highest monomial of $\mathcal{W}_{k,a}^{(i)}$ and $\lambda = k\Lambda_i - \sum_{j\in I_n}\alpha_j$. Then
$$\text{dim}((\mathcal{W}_{k,a}^{(i)})_{\lambda}) = n+1.$$
Indeed from Theorem \ref{nexp} this weight space corresponds to the $n$ monomials
$$MA_{i,q^{2k-1}}^{-1}(A_{i+1,aq^{2k}}^{-1}\cdots A_{i+j,aq^{2k+j-1}}^{-1})(A_{i-1,aq^{2k}}^{-1}\cdots A_{i-j',aq^{2k+j'-1}}^{-1})$$
where $j+j' = n$, $0\leq j,j'$. The monomials for $j-1 = j' = (n-1)/2$ and $j'-1 = j = (n-1)$ are the same. But this monomial occurs with multiplicity $2$. The simple module $L$ of $\U_q^h$ of highest weight $k\Lambda_i$ satisfies $\text{dim}(L_\lambda) = n$.

Now suppose that $n$ is even. Let $\lambda = k\Lambda_i - 2\sum_{j\in I_n}\alpha_j$. Then
$$\text{dim}((\mathcal{W}_{k,a}^{(i)})_{\lambda})
= 2n + 2 + (n+1)(n+2)/2.$$
Indeed from Theorem \ref{nexp} this weight space corresponds to the $2n + 1 + (n+1)(n+2)/2$ monomials 
$$MA_{i,q^{2k-1}}^{-1}(A_{i+1,aq^{2k}}^{-1}\cdots A_{i+j,aq^{2k+j-1}}^{-1})(A_{i-1,aq^{2k}}^{-1}\cdots A_{i-j',aq^{2k+j'-1}}^{-1})$$
$$\times A_{i,q^{2k-3}}^{-\epsilon}(A_{i+1,aq^{2k-2}}^{-1}\cdots A_{i+l,aq^{2k+l-3}}^{-1})(A_{i-1,aq^{2k}}^{-1}\cdots A_{i-j',aq^{2k+j'-1}}^{-1})$$
where $\epsilon\in\{0,1\}$, $0\leq l\leq j$, $0\leq l'\leq j'$, $j+j'+\epsilon + l + l' = 2n + 1$, $(\epsilon = 0\Rightarrow l = l' = 0)$, $(\epsilon = 1\Rightarrow l + j' = l + j = n)$. 

\noindent The term $2n + 1$ corresponds to the case $\epsilon = 0$, and the term $(n+1)(n+2)/2$ to the case $\epsilon = 1$ as we have to count the number of couples $(j,j')$ satisfying $0\leq j,j'\leq n$ and $j+j'\geq n$ 

\noindent The monomials for $\epsilon = 0$, $j'- 1 = j = n$ and $j - 1 = j' = n$ are the same. But this monomial occurs with multiplicity $2$.

\noindent The simple module $L$ of $\U_q^h$ of highest weight $k\Lambda_i$ satisfies $\text{dim}(L_\lambda) = 2n+1+(n+1)(n+2)/2$.
\qed

\begin{rem} In particular the module $\mathcal{W}_{k,a}^{(i)}$ can not be obtained by evaluation
from a representation of $\U_q^h$ (compare with Remark \ref{eval}).\end{rem}

\noindent For example consider $\U_q(sl_4^{tor})$ and $\mathcal{W}_{1,1}^{(0)}$.
The first terms of the $q$-character as computed in \cite{Selecta} are $Y_{0,1}
+ Y_{0,q^2}^{-1}Y_{1,q}Y_{3,q}
+ Y_{3,q^3}^{-1}Y_{1,q}Y_{2,q^2} + Y_{1,q^3}^{-1}Y_{2,q^2}Y_{3,q}
+ Y_{2,q^4}^{-1}Y_{1,q}Y_{1,q^3}
+ Y_{1,q^3}^{-1}Y_{2,q^2}^2Y_{3,q^3}^{-1}Y_{0,q^2}
+ Y_{2,q^4}^{-1}Y_{3,q}Y_{3,q^3}
+ Y_{1,q}Y_{1,q^5}^{-1}Y_{0,q^4} + Y_{0,q^4}^{-1}Y_{2,q^2}^2
+ 2\times Y_{2,q^2}Y_{2,q^4}^{-1}Y_{0,q^2}
+Y_{3,q^5}^{-1}Y_{3,q}Y_{0,q^4}+\cdots$. The $4$ monomials (with multiplicity) corresponding to the weight $\lambda = \Lambda_0 - \alpha_0 - \alpha_1 - \alpha_2 - \alpha_3$ are $Y_{1,q}Y_{1,q^5}^{-1}Y_{0,q^4}+ 2\times Y_{2,q^2}Y_{2,q^4}^{-1}Y_{0,q^2}
+Y_{3,q^5}^{-1}Y_{3,q}Y_{0,q^4}$.

\section{Positivity}\label{pos}

In this section we address several natural questions related to the results of this note.
In particular we write a general conjecture extending Theorem \ref{nexp} (Conjecture \ref{cpos}).

The algebra $sl_\infty$ has weight lattice $P$ and we can define its category $\mathcal{O}_{\text{int}}$
as for $\U_q(sl_\infty)$. There is a connection between the representation theory
of $sl_\infty$ and of the affine Lie algebras $\hat{sl}_{n+1}$.

We have a natural projection map $\phi_n : P\rightarrow P_n$ where $P_n = \sum_{i\in I_n}\ZZ \Lambda_i$
is the weight lattice of $\hat{sl}_{n+1}$ : for $i\in\ZZ$ we set $\phi_n(\Lambda_i) = \Lambda_{[i]}$.
The map $\phi_n$ is naturally extended to characters of objects in the category $\mathcal{O}_{\text{int}}$
of $sl_\infty$. But we can not get direct character formulae as we did in Section \ref{torapp}.

For example, consider the integrable simple representation $L$ of $\hat{sl}_4$
of highest weight $\Lambda_0$ (in the standard sense, that is with respect
to Chevalley generators). For $\lambda = \Lambda_0 - \alpha_0 - \alpha_{-1} - \alpha_1 - \alpha_2$,
we have $\text{dim}(L_\lambda) = 3$.
But $\phi_n^{-1}(\lambda)$ contains the following weights : 
$\Lambda_0 - \alpha_0 - \alpha_{-1} - \alpha_{-2} - \alpha_{-3}$,
$\Lambda_0 - \alpha_0 - \alpha_{-1} - \alpha_{-2} - \alpha_1$,
$\Lambda_0 - \alpha_0 - \alpha_1 - \alpha_2 - \alpha_{-1}$,
$\Lambda_0 - \alpha_0 - \alpha_1 - \alpha_2 - \alpha_3$. 
Each of them corresponds to a non trivial weight space in the simple integrable 
representation $L'$ of $sl_\infty$ of highest weight $\Lambda_0$. 
So $\phi_n(\chi(L'))\neq \chi_n(L)$ where $\chi_n$ is the character map for $\hat{sl}_{n+1}$.

However, we clearly have $\phi_n(\text{Im}(\chi))\subset \text{Im}(\chi_n)$
by using the characterization with the Weyl group action. So for $V$ simple in 
$\mathcal{O}_{\text{int}}$ for $sl_\infty$, we have
a decomposition
$$\phi_n(\chi(V)) = \sum_W n_W \chi_n(W)$$
where the sum is over the isomorphism classes $W$ of simple integrable representations
of $\hat{sl}_{n+1}$ in the category $\mathcal{O}$. The $n_W\in\ZZ$ and $\phi_n(\chi(V))$ is the character of a representation of $\hat{sl}_{n+1}$, which is a priori virtual. But we have the following positivity result.

\begin{prop}\label{posres} $\phi_n(\chi(V))$ is the character of an actual representation of $\hat{sl}_{n+1}$,
that is, for any $W$, we have $n_W\geq 0$.\end{prop}

\demo Consider the formal infinite sum of elements of $sl_\infty$ for $i\in I_n$
$$\hat{x}_i^\pm = \sum_{r\in\ZZ} x_{i + (n+1) r}^\pm
\text{ , }\hat{h}_i = \sum_{r\in\ZZ} h_{i + (n+1) r}.$$
Then consider $V$ simple in $\mathcal{O}_{\text{int}}$ for $sl_\infty$.
Then $\hat{x}_i^\pm, \hat{h}_i$ make sense as operators on $V$ : for weight reason,
for each $v\in V$, only a finite number of the terms of the infinite sum have a
non zero action on $v$. Then it easy to check that they satisfy the relations of $\hat{sl}_{n+1}$.
The less obvious relations to be checked are the Serre relations.
First suppose that $n\geq 2$. For $i\in \ZZ$, we get
$$[\hat{x}_{[i]+1}^+,[\hat{x}_{[i]+1}^+,\hat{x}_{[i]}^+]]
= \sum_{r,r',r''\in\ZZ} [x_{i + 1 + (n+1) r}^+, [x_{i + 1 + (n+1)r'}^+ ,[x_{i + (n + 1)r''}^+]]]$$
$$= \sum_{r\in\ZZ} [x_{i + 1 + (n+1) r}^+, [x_{i + 1 + (n+1)r}^+ ,[x_{i + (n + 1)r}^+]]]
= 0.$$
Now suppose that $n = 1$. We get
$$[\hat{x}^+_{[i]+1},[\hat{x}^+_{[i]+1},[\hat{x}^+_{[i]+1},\hat{x}^+_{[i]}]]]
= \sum_{r,r',r'',r'''\in\ZZ} [x^+_{i+1+2r},[x^+_{i + 1 + 2 r'}, [x^+_{i + 1 + 2r''} ,x^+_{i + 2r'''}]]]$$
$$= \sum_{r\in\ZZ,\epsilon,\epsilon', \epsilon'' \in\{1,-1\}} [x^+_{i+\epsilon'' +2r},[x^+_{i + \epsilon' + 2r}, [x^+_{i + \epsilon + 2r} ,x^+_{i + 2r}]]]$$
$$= \sum_{r\in\ZZ,\epsilon, \epsilon'' \in\{1,-1\}} [x^+_{i+\epsilon'' +2r},[x^+_{i - \epsilon + 2r}, [x^+_{i + \epsilon + 2r} ,x^+_{i + 2r}]]]$$
$$=
- \sum_{r\in\ZZ,\epsilon, \epsilon'' \in\{1,-1\}} [x^+_{i-\epsilon +2r},[x^+_{i + \epsilon'' + 2r}, [x^+_{i + \epsilon + 2r} ,x^+_{i + 2r}]]]$$
$$=
- \sum_{r\in\ZZ,\epsilon\in\{1,-1\}} [x^+_{i+\epsilon +2r},[x^+_{i - \epsilon + 2r}, [x^+_{i - \epsilon + 2r} ,x^+_{i + 2r}]]]
= 0.$$
So we get an action of $\hat{sl}_{n+1}$ on $V$, and the character of this representation
is precisely $\phi_n(\chi(V))$.
\qed

\begin{rem}\label{knclass} Analogs of the operators $\hat{x}^\pm_i, \hat{h}_i$ are considered in
\cite[Proposition 3.5]{jmmo} for Fock spaces. Moreover, it was pointed out to the
author by Edward Frenkel that there is a simple reason for these operators to
satisfy the relations of $\hat{sl}_{n+1}$. Let $V$ be the natural $n$-dimensional
representation of $sl_n$. Then $sl_n[t^{\pm 1}]$ acts on $V[t^{\pm 1}]$. By choosing a basis of weight vectors, $V[t^{\pm 1}]$ can be identified with $\CC[t^{\pm 1}]$ and
we get a Lie algebra morphism $\Psi : sl_n[t^{\pm 1}]\rightarrow sl(\CC[t^{\pm 1}])$. Here
$sl(\CC[t^{\pm 1}])$ can be seen as a completion of $sl_\infty$. This map $\Psi$ is precisely given by the explicit
formulae considered in the proof of the Proposition \ref{posres}.
\end{rem}

This result leads to a natural positivity conjecture in the context of the results of this note. Let $V$ be a simple representations in $\mathcal{O}_{\text{int}}$ for $\U_q(\hat{sl}_\infty)$.
Then we can prove as in Theorem \ref{nexp} that $\phi_n(\chi_q(V))$ is the $q$-character of a representation
of $\U_q(sl_{n+1}^{tor})$, which is a priori virtual.

\begin{conj}\label{cpos} $\phi_n(\chi_q(V))$ is the $q$-character of an actual representation
of $\U_q(sl_{n+1}^{tor})$.
\end{conj}

A priori this Conjecture can not be proved as Proposition \ref{posres} as there is no
obvious action of $\U_q(sl_{n+1}^{tor})$ on $V$, as far the author can
see.
However, Theorem \ref{nexp} gives a proof of this Conjecture 
for Kirillov-Reshetikhin modules.

In general, the corresponding (conjectural) representation 
of $\U_q(sl_{n+1}^{tor})$ will not be simple.
For example consider $V$ of highest monomial $Y_{0,1}Y_{2,q^4}$. Then $\chi_q(V)$
contains the monomial $Y_{0,1}Y_{2,q^4}A_{0,q}^{-1}A_{-1,q^2}^{-1}A_{-2,q^3}^{-1}$.
This follows from Theorem \ref{limit} as for the sub quantum affine algebra
we get a minimal affinization whose $q$-character is known (see references 
in the introduction of \cite{h4}).
Its image by $\phi_n$ is $Y_{1,q}Y_{1,q^3}$ which does not occur in the $q$-character
of the simple $\U_q(sl_4^{tor})$-module $V_1$ of highest monomial $Y_{0,1}Y_{2,q^4}$.
This is a particular case of the elimination theorem \cite[Theorem 5.1]{AnnENS}.
In fact the dominant monomials in $\phi_n(\chi_q(V))$ are $Y_{0,1}Y_{2,q^4}$, $Y_{2,q^2}$
and $Y_{1,q}Y_{1,q^3}$ and they occur with multiplicity $1$. As the first two ones occur in
the $q$-character of $V_1$ with multiplicity $1$, we get
$$\phi_n(\chi_q(V)) = \chi_q(V_1) + \chi_q(V_2)$$
where $V_2$ is the simple $\U_q(sl_4^{tor})$-module $V_1$ of highest monomial $Y_{0,1}Y_{2,q^4}$.
In particular the positivity Conjecture \ref{cpos} holds in this case.

\end{document}